\documentclass[11pt]{amsart}
\usepackage[utf8]{inputenc}
\usepackage[T1]{fontenc}
\usepackage{amsmath,amssymb,amsthm,mathrsfs,enumitem}
\usepackage[margin=2.7cm]{geometry}
\usepackage{xcolor}
\usepackage[colorlinks=true,linkcolor=blue!50!black,citecolor=blue!50!black,urlcolor=blue!50!black]{hyperref}

\newtheorem{theorem}{Theorem}[section]
\newtheorem{maintheorem}{Theorem}

\newtheorem*{thmBprime}{Theorem B$'$}
\newtheorem{proposition}[theorem]{Proposition}
\newtheorem{lemma}[theorem]{Lemma}
\newtheorem{corollary}[theorem]{Corollary}
\theoremstyle{definition}

\newtheorem{remark}[theorem]{Remark}
\numberwithin{equation}{section}

\newcommand{\Q}{\mathbb{Q}}
\newcommand{\Z}{\mathbb{Z}}
\newcommand{\F}{\mathbb{F}}
\newcommand{\T}{\mathbb{T}}
\newcommand{\cO}{\mathcal{O}}
\newcommand{\GL}{\mathrm{GL}}
\newcommand{\SL}{\mathrm{SL}}
\newcommand{\Ind}{\mathrm{Ind}}
\newcommand{\Gal}{\mathrm{Gal}}
\newcommand{\Ord}{\mathrm{Ord}}
\newcommand{\Spec}{\mathrm{Spec}\,}
\newcommand{\Hom}{\mathrm{Hom}}
\newcommand{\Art}{\mathrm{Art}}
\newcommand{\HT}{\mathrm{HT}}
\newcommand{\ps}{\mathrm{ps}}
\newcommand{\cl}{\mathrm{cl}}
\newcommand{\nc}{\mathrm{nc}}
\newcommand{\ad}{\mathrm{ad}}
\newcommand{\rhobar}{\bar\rho}

\newcommand{\So}{\Sigma_{\mathrm{o}}}
\newcommand{\Sn}{\Sigma_{\mathrm{n}}}

\title[Some instances of modularity in the mixed case]{Some instances of modularity for residually reducible representations in the mixed case}
\author{Luis V. Dieulefait}

\address{Departament de Matem\`atiques i Inform\`atica,
  Universitat de \break Barcelona, Gran Via de les Corts Catalanes,
  585. 08007 Barcelona}
\email{ldieulefait@ub.edu}
\begin{document}

\keywords{Fontaine-Mazur conjecture}
\subjclass[2010]{11F80, 11F33}

\begin{abstract}
Let $F$ be a totally real field, abelian over $\Q$, in which an odd prime $p$ splits completely, and let $\rho:G_F\to\GL_2(\overline\Q_p)$ be an odd irreducible representation which is de Rham with distinct Hodge--Tate weights at the places above $p$ and whose residual representation is reducible. We consider the \emph{mixed} case of the Fontaine--Mazur conjecture, in which $\rho$ is reducible at some places above $p$ and irreducible at the others. This case is not covered by the work of Pan and Zhang. We prove some instances of modularity in this situation. Our main result is that $\rho$ is modular if it is non-split at every place above $p$ where it is reducible, by combining Zhang's pro-modularity theorem with a classicality criterion based on Emerton's ordinary parts. Together with a theorem of Zhao on local indecomposability, this gives the modularity of the $\lambda$-adic representations attached to abelian varieties of $\GL_2$-type in the mixed case, and, using a result of Snowden, of potentially modular representations of parallel weight two under a mild additional condition. It also gives the modularity of representations whose local type at every place above $p$ is of discrete series type, and, through a member of the compatible system with interchanged Hodge--Tate weights, of potentially modular representations of weight $(2,k)$ over real quadratic fields with Steinberg type at one place above $p$, under conditions on that member.
\end{abstract}

\maketitle
\tableofcontents

%======================================================================
\section{Introduction}

\subsection{The mixed case}
Let $F$ be a totally real number field, $p$ a prime and $\rho:G_F\to\GL_2(\overline\Q_p)$ a continuous irreducible representation, unramified outside finitely many primes. The Fontaine--Mazur conjecture predicts that if $\rho$ is odd and $\rho|_{G_{F_v}}$ is de Rham with distinct Hodge--Tate weights for every $v\mid p$, then $\rho$ is \emph{modular}, i.e.\ it arises from a regular algebraic cuspidal automorphic representation of $\GL_2(\mathbb A_F)$.

When the residual representation $\rhobar$ is reducible, the first general results are due to Skinner and Wiles \cite{SW99}, under the assumption that $\rho$ is ordinary at all places above $p$. Pan \cite{Pan22} introduced a method based on the completed cohomology of definite quaternion algebras and Paškūnas' theory of blocks \cite{Pas13}, which treats the case where $F$ is abelian over $\Q$, $p$ splits completely in $F$ and $\rho|_{G_{F_v}}$ is irreducible for every $v\mid p$. Zhang \cite{Zha25} removed the remaining exceptional residual configuration at $p=3$, using the results of Paškūnas--Tung \cite{PT21}. The \emph{mixed} case, in which $\rho|_{G_{F_v}}$ is reducible at some places $v\mid p$ and irreducible at the others, is not covered by these works. In this paper we prove some instances of modularity in the mixed case.

\subsection{Setting}\label{sec:setting}
Throughout, $p$ is an odd prime and $F$ is a totally real field, abelian over $\Q$, in which $p$ splits completely, so that $F_v=\Q_p$ for every $v\mid p$. Let $E/\Q_p$ be a finite extension with ring of integers $\cO$ and residue field $\F$, and let
\[
\rho:G_F\longrightarrow\GL_2(E)
\]
be a continuous, absolutely irreducible, odd representation, unramified outside finitely many places. We write $\rho_v=\rho|_{G_{F_v}}$, $\Sigma_p$ for the set of places of $F$ above $p$, and we fix a partition $\Sigma_p=\Sn\sqcup\So$. We consider the following conditions.
\begin{itemize}
\item[(N)] For $v\in\Sn$, $\rho_v$ is irreducible and de Rham with distinct Hodge--Tate weights.
\item[(O)] For $v\in\So$, $\rho_v$ is reducible and de Rham with distinct Hodge--Tate weights.
\item[(R)] $\rhobar^{ss}=\bar\chi_1\oplus\bar\chi_2$ for characters $\bar\chi_i:G_F\to\F^\times$, the character $\bar\chi=\bar\chi_1\bar\chi_2^{-1}$ extends to a character of $G_\Q$, and $\bar\chi|_{G_{F_v}}\neq1$ for every $v\in\Sigma_p$.
\end{itemize}
Here and below we freely extend scalars: a local representation $\rho_v$ which is irreducible but not absolutely irreducible becomes reducible over a finite extension of $E$, so after enlarging $E$ each $\rho_v$ is either reducible or absolutely irreducible, and \emph{irreducible} means absolutely irreducible. Condition (R) does not depend on the order of $\bar\chi_1,\bar\chi_2$. Under (O), the two Jordan--Hölder characters of $\rho_v$ have distinct Hodge--Tate numbers. By Bloch--Kato theory (Lemma~\ref{lem:BK}), the one with smaller Hodge--Tate number, in the convention $\HT(\varepsilon)=-1$, is always a subcharacter of $\rho_v$, and it is the only one if $\rho_v$ is non-split. We call it the \emph{classical character} of $\rho_v$; it is the subcharacter which occurs for ordinary modular forms.

\subsection{The non-split case}
Our first result treats the case where $\rho$ is non-split at all the places where it is reducible.

\begin{maintheorem}\label{thm:A}
Assume \textup{(N)}, \textup{(O)} and \textup{(R)}, and that $\rho_v$ is non-split for every $v\in\So$. Then $\rho$ is modular.
\end{maintheorem}

When $\So=\emptyset$ this is the theorem of Pan and Zhang, and when $\Sn=\emptyset$ it is covered by the work of Skinner--Wiles in parallel weight, and in general by \cite[Theorem~5.1.2]{Pan22}; the new content is the mixed case. No condition is imposed on the Hodge--Tate weights beyond (N) and (O); in particular they need not be parallel. The hypothesis that $\rho_v$ is non-split for $v\in\So$ is the natural one for our classicality criterion: if $\rho_v$ is split, the method only shows that \emph{one} of its two subcharacters occurs in the relevant Hida family, and the non-classical one gives no classical vectors.

We should be clear about what is new in the proof of Theorem~\ref{thm:A}. The pro-modularity input is \cite[Theorem~5.5.1]{Zha25} (for $p\geq5$, the corresponding result of Pan \cite{Pan22}), applied verbatim: its proof does not use any local condition at $p$, so it applies in the mixed case, although this is not stated there. The ingredients of the classicality step are standard as well: Emerton's ordinary parts and his adjunction formula \cite{Eme10a,Eme10b}, the local--global compatibility of Pan and Zhang for Paškūnas blocks \cite{Pan22,Zha25}, the results of Breuil--Emerton on locally algebraic vectors in ordinary principal series \cite{BE10}, and those of Colmez--Dospinescu--Paškūnas at the non-ordinary places \cite{CDP14}. The only new point is to combine them, taking ordinary parts only at the places of $\So$; this realizes, at the level of points, the mixture of completed cohomology and Hida families suggested by Pan \cite[Remark~8.0.3]{Pan22}, and we do not claim any originality for it. Moreover, the non-split hypothesis avoids precisely the genuine difficulty of the mixed case, namely showing that the classical character occurs. The main interest of the paper lies instead in the applications of Theorem~\ref{thm:A} given in Sections~\ref{sec:appA} and~\ref{sec:AV}, all the more so since they will be used to prove new modularity results that do not depend on any extension, not yet available, of Pan's theorem to the mixed case (see Remark~\ref{rem:applications}).

\begin{remark}\label{rem:GreenbergA}
Under the hypothesis of potential modularity, which is available in the applications we have in mind (Remark~\ref{rem:applications}), Greenberg's question is exactly what is needed to deduce from Theorem~\ref{thm:A} the modularity of $\rho$ in the mixed case in general, that is, without assuming that $\rho_v$ is non-split for $v\in\So$. Indeed, if $\rho$ is induced from a character $\theta$ of $G_K$ for a quadratic extension $K/F$, then it is easy to see that it is modular. Assume that $\rho$ is not induced, and that $\rho|_{G_W}$ is modular for a totally real Galois extension $W/F$. Then it is not hard to see that $\rho|_{G_W}$ is not induced either, so it corresponds to a Hilbert newform with no complex multiplication. A question of Greenberg (see \cite{GV04}) predicts that the Galois representation attached to such a Hilbert eigenform is non-split at every place above $p$ at which it is ordinary. Applied to $\rho|_{G_W}$ at a place $w$ above $v\in\So$, it would give that $\rho_v$ is non-split, so that Theorem~\ref{thm:A} applies. Thus, under potential modularity, a positive answer to Greenberg's question for Hilbert newforms would, together with Theorem~\ref{thm:A}, give the modularity of $\rho$ in the mixed case in general. Potential modularity is essential here: without it, the Fontaine--Mazur conjecture together with Greenberg's question would also predict that $\rho_v$ is never split for $v\in\So$ outside the CM case, but modularity is precisely what we want to prove, so we cannot assume it, whereas potential modularity is a hypothesis that is available in applications. This explains why several applications of Theorem~\ref{thm:A} can be found (see Sections~\ref{sec:appA} and~\ref{sec:AV}), and the more cases of Greenberg's question are settled, the more general this result becomes.
\end{remark}

\subsection{Three applications of Theorem~\ref{thm:A}}
Theorem~\ref{thm:A} also applies whenever the local types at the places above $p$ force each $\rho_v$ to be either irreducible or reducible and non-split. Say that $\rho_v$ is of discrete series type if its Frobenius-semisimple Weil--Deligne representation is indecomposable, that is, either irreducible (supercuspidal type) or a twist of the Steinberg one (Steinberg type, i.e.\ $\rho_v$ is a twist of a semistable non-crystalline representation). In the supercuspidal case $\rho_v$ is irreducible, and in the Steinberg case $\rho_v$ is reducible and non-split if its Hodge--Tate numbers differ by $1$, and irreducible otherwise (Lemma~\ref{lem:St}). Hence $\rho$ is modular if it satisfies (R) and, for every $v\in\Sigma_p$, $\rho_v$ has distinct Hodge--Tate numbers and is of discrete series type (Corollary~\ref{cor:St}); over a real quadratic field it suffices that the Hodge--Tate numbers differ by $1$ at both places above $p$ and that $\rho$ is of Steinberg type at one of them (Corollary~\ref{cor:St2}). In non-parallel weight $(2,k)$ over a real quadratic field, with Steinberg type at one of the two places above $p$ and under potential modularity, we pass to a conjugate member $\rho'$ of the compatible system containing $\rho$, in which the Hodge--Tate numbers at the two places are interchanged. Theorem~\ref{thm:2k} reduces the modularity of $\rho$ either to a residually reducible case covered by Theorem~\ref{thm:A} and \cite[Theorem~5.1.2]{Pan22}, or to a residually irreducible case covered by the modularity lifting theorems of Kisin, Hu--Tan and Tung.

\subsection{Abelian varieties of $\GL_2$-type}
Theorem~\ref{thm:A} applies to abelian varieties. Let $A/F$ be an abelian variety of $\GL_2$-type: there is a number field $K\subset\mathrm{End}_F(A)\otimes\Q$ with $[K:\Q]=\dim A$. For a prime $\lambda$ of $K$ above $p$, the $\lambda$-adic Tate module of $A$ gives a representation
\[
\rho_{A,\lambda}:G_F\longrightarrow\GL_2(K_\lambda).
\]

\begin{maintheorem}\label{thm:B}
Assume that $\rho_{A,\lambda}$ satisfies \textup{(R)}. Then $\rho_{A,\lambda}$ is modular: it arises from a cuspidal Hilbert eigenform of parallel weight $2$ over $F$.
\end{maintheorem}

No irreducibility assumption is needed: since $F$ is totally real, $\rho_{A,\lambda}$ is always absolutely irreducible, by a well-known argument of Ribet based on purity (Lemma~\ref{lem:AVirr}). No assumption is made on the reduction of $A$ at the places above $p$. The places $v\mid p$ where $\rho_{A,\lambda}$ is reducible (after extending scalars) form the set $\So$, and the others the set $\Sn$. If $A$ has complex multiplication over $\overline\Q$, modularity follows from the theory of complex multiplication and the automorphic induction of Hecke characters. Otherwise, a theorem of Zhao \cite[Theorem~5.2]{Zhao14}, which extends a theorem of Hida \cite{Hid13}, shows that $\rho_{A,\lambda}|_{G_{F_v}}$ is non-split for every $v\in\So$, and Theorem~\ref{thm:A} applies.

As a consequence, we obtain in Section~\ref{sec:Bprime} the modularity of potentially modular representations of parallel weight $2$ in the mixed case (Theorem~B$'$), under the condition that either $[F:\Q]$ is odd or the Weil--Deligne representation of $\rho$ at some finite place (possibly above $p$) is indecomposable, i.e.\ corresponds to a discrete series representation. Such a representation lies in a strongly compatible system, by the descent construction of \cite{Die04} and of the proof of \cite[Theorem~5.5.1]{BLGGT14}, and under this condition a result of Snowden \cite{Sno09}, applied to a member of the system with large residual image, shows that it comes from an abelian variety of $\GL_2$-type over $F$, so that Theorem~\ref{thm:B} applies.

\begin{remark}\label{rem:applications}
Although the results of the present paper do not settle the mixed case of Pan's theorem in general, they are sufficient for significant applications to new modularity results. In particular, two of the results of this paper, Theorem~B$'$ and Corollary~\ref{cor:St}, are already used in \cite{DP26}. In such applications, and more generally in the propagation of modularity along safe chains of congruences, the hypothesis of potential modularity in Theorems~B$'$ and~\ref{thm:2k} is available whenever the construction already supplies a potentially modular compatible system, as was the case in the proofs of Serre's modularity conjecture over $\Q$ by Khare--Wintenberger \cite{KW09a,KW09b} and by the author \cite{Die12b}. In inductive arguments of this kind, the residual hypotheses of Theorem~\ref{thm:2k}\textup{(c)} on the entangled representation $\rho'$ are also provided by the induction, since the invariants of $\rho'$ are as small as those of $\rho$ (see Remark~\ref{rem:entangled}).
\end{remark}

\subsection{Strategy}
After a solvable base change $F_1/F$ provided by Zhang, the pro-modularity theorem of \cite{Zha25} shows that $\rho|_{G_{F_1}}$ comes from a maximal ideal $\mathfrak p$ of a big Hecke algebra acting on the completed cohomology $S$ of a definite quaternion algebra. Its proof does not use any local condition at $p$. The task is to show that this point is classical.

\medskip\noindent\textbf{Ordinary parts and classicality.} We apply Emerton's functor of ordinary parts \cite{Eme10a,Eme10b} to $S$, with respect to the Borel subgroup at the places of $F_1$ above $\So$ only. The resulting space $S^{\mathrm o}$ is completed cohomology at the places above $\Sn$ and a Hida family at the places above $\So$, and its Hecke algebra $\T^{\mathrm o}$ carries characters $\psi_{v,1}$ given by the diamond operators. We show that the pseudo-representation over $\T^{\mathrm o}$ is $\psi_{v,1}$-ordinary (Lemma~\ref{lem:galois}); since the central character need not be locally constant at $p$, this is first proved after a twist by a continuous character and then transported back. Starting from an irreducible subrepresentation of the $\mathfrak p$-eigenspace of $S$, Emerton's adjunction produces ordinary eigenvectors on which the torus acts through a character whose first component is \emph{a} subcharacter of $\rho|_{G_{F_{1,v}}}$ (Proposition~\ref{prop:LG}). If it is the classical character, the corresponding principal series and the non-ordinary local factors have dense locally algebraic vectors, and the adjunction map produces classical forms (Theorem~\ref{thm:class}). If $\rho_v$ is non-split, its only subcharacter is the classical one; this proves Theorem~\ref{thm:A}.

\medskip\noindent\textbf{Abelian varieties.} For Theorem~\ref{thm:B}, the main additional input is the local indecomposability theorem of Hida and Zhao, which is a statement about abelian varieties and does not use modularity; the CM case is classical. Theorem~B$'$ follows from Theorem~\ref{thm:B}, since by a result of Snowden, applied to a suitable member of the compatible system containing $\rho$, the representations it concerns come from abelian varieties of $\GL_2$-type.

\subsection{Organisation}
Section~\ref{sec:local} contains local results on reducible representations of $G_{\Q_p}$. In Section~\ref{sec:aut} we recall completed cohomology and Paškūnas' blocks, and develop the ordinary parts at the ordinary places. Section~\ref{sec:class} contains the local--global compatibility at ordinary points, the classicality criterion, the pro-modularity input, and the proof of Theorem~\ref{thm:A}. Section~\ref{sec:appA} contains three applications of Theorem~\ref{thm:A}: to representations of discrete series type at all places above $p$ (Corollary~\ref{cor:St}), to representations of parallel weight $2$ over real quadratic fields with Steinberg type at one place above $p$ (Corollary~\ref{cor:St2}), and, through entangled pairs, to representations of weight $(2,k)$ (Theorem~\ref{thm:2k}). Section~\ref{sec:AV} is devoted to abelian varieties: it contains the proof of Theorem~\ref{thm:B} and of its consequence Theorem~B$'$ for potentially modular representations of parallel weight $2$.

\subsection{Acknowledgments}
The author acknowledges the use of Claude (Anthropic), an AI assistant, during the preparation of this work: it was used to check the relevant parts of the literature, to work out and verify technical details of the arguments and to help draft the manuscript. The author has checked all the mathematical content and takes full responsibility for it.

\subsection*{Notation}
For a field $k$ we write $G_k$ for its absolute Galois group, and for a number field $K$ and a finite set $S$ of places we write $G_{K,S}$ for the Galois group of the maximal extension unramified outside $S$. We always include the infinite places in $S$, and we write $G_v=G_{K_v}$ for the decomposition group at a place $v$. We use the convention $\HT(\varepsilon)=-1$ for Hodge--Tate numbers, as in \cite{Pan22,Zha25}. Characters of $G_{\Q_p}$ are identified with characters of $\Q_p^\times$ via the local Artin map $\Art$, normalized as in \cite{Pan22,Zha25}, so that $\varepsilon\circ\Art(x)=x|x|_p$; in particular $\varepsilon\circ\Art(u)=u$ for $u\in\Z_p^\times$. We write $B$ (resp.\ $\bar B$) for the upper (resp.\ lower) Borel subgroup of $\GL_2$, $N$ for the unipotent radical of $B$ and $T$ for the diagonal torus. For a representation $V$ we write $\ad^0V$ for the trace zero endomorphisms, and $V(1)=V\otimes\varepsilon$. For a totally real field $K$, a regular algebraic cuspidal automorphic representation $\pi$ of $\GL_2(\mathbb A_K)$, a prime $\ell$ and an isomorphism $\iota:\overline\Q_\ell\simeq\mathbb C$, we write $r_{\pi,\iota}:G_K\to\GL_2(\overline\Q_\ell)$ for the Galois representation attached to $\pi$. A representation $\rho$ of $G_K$ is \emph{modular} if $\rho\otimes\overline\Q_\ell\simeq r_{\pi,\iota}$ for some such $\pi$ and $\iota$. We enlarge $E$ when necessary.

%======================================================================
\section{Local preliminaries}\label{sec:local}

In this section $G=G_{\Q_p}$.

\subsection{Reducible de Rham representations}

\begin{lemma}\label{lem:BK}
Let $\psi_1,\psi_2:G\to E^\times$ be de Rham characters with $\HT(\psi_1)\neq\HT(\psi_2)$.
\begin{enumerate}
\item If $\HT(\psi_1)<\HT(\psi_2)$, every extension of $\psi_2$ by $\psi_1$ is de Rham.
\item If $\HT(\psi_1)>\HT(\psi_2)$, every de Rham extension of $\psi_2$ by $\psi_1$ is split.
\end{enumerate}
\end{lemma}

\begin{proof}
Let $\eta=\psi_1\psi_2^{-1}$. Extension classes are given by $H^1(G,E(\eta))$, and de Rham ones by the Bloch--Kato subspace $H^1_g$. By the dimension formulas of \cite[Corollary~3.8.4]{BK90} and local duality, $H^1_g=H^1$ when $\HT(\eta)<0$, while $H^1_g=0$ when $\HT(\eta)>0$.
\end{proof}

Let $\rho_v:G\to\GL_2(E)$ be reducible and de Rham with distinct Hodge--Tate weights, with Jordan--Hölder characters $\psi^{\cl},\psi^{\nc}$ such that $\HT(\psi^{\cl})<\HT(\psi^{\nc})$. By Lemma~\ref{lem:BK}, $\psi^{\cl}$ is a subcharacter of $\rho_v$, and it is the only one if $\rho_v$ is non-split. We call it the \emph{classical character}. It is the subcharacter occurring for ordinary modular forms.

\begin{lemma}\label{lem:la}
Let $\mu=\mu_1\otimes\mu_2$ be a unitary character of $T(\Q_p)$ with $\mu_i|_{\Z_p^\times}=x^{a_i}\theta_i$, $a_i\in\Z$ and $\theta_i$ of finite order. Then $(\Ind_{\bar B}^{\GL_2(\Q_p)}\mu)_{\mathrm{cont}}$ has nonzero locally algebraic vectors if and only if $a_1\geq a_2$, in which case they are dense. In particular, if $\psi^{\cl},\psi^{\nc}$ are as above and $\mu^{\cl}=\psi^{\cl}\otimes\psi^{\nc}\varepsilon$, then $(\Ind_{\bar B}\mu^{\cl})_{\mathrm{cont}}$ has dense locally algebraic vectors.
\end{lemma}

\begin{proof}
Algebraic functions $f$ on $\GL_2$ with $f(\bar bg)=\mu_{\mathrm{alg}}(\bar b)f(g)$ are, up to ${\det}^{a_2}$, homogeneous polynomials of degree $a_1-a_2$ in the first row, and there are none if $a_1<a_2$. If $a_1\geq a_2$, the locally algebraic vectors are $W\otimes(\Ind_{\bar B}\mu_{\mathrm{sm}})^{\mathrm{sm}}$, with $W=\mathrm{Sym}^{a_1-a_2}\otimes{\det}^{a_2}$ and $\mu_{\mathrm{sm}}=\mu\mu_{\mathrm{alg}}^{-1}$. For density, view $(\Ind_{\bar B}\mu)_{\mathrm{cont}}$ as the space of continuous sections of a line bundle on $\bar B\backslash\GL_2(\Q_p)=\mathbf P^1(\Q_p)$. Cover $\mathbf P^1(\Q_p)$ by finitely many disjoint compact open subsets $V_i$ on each of which there is a nowhere vanishing section $s_i$ coming from $W\otimes(\Ind_{\bar B}\mu_{\mathrm{sm}})^{\mathrm{sm}}$. A continuous section is $\sum_if_is_i$ with $f_i$ continuous on $V_i$, and approximating the $f_i$ uniformly by locally constant functions gives locally algebraic approximations (cf.\ \cite[Proposition~2.2.1]{BE10} for $a_1>a_2$). With our conventions $a_1=-\HT(\psi^{\cl})$ and $a_2=1-\HT(\psi^{\nc})$, so $a_1-a_2=\HT(\psi^{\nc})-\HT(\psi^{\cl})-1\geq0$.
\end{proof}

\subsection{Ordinary lines}

The following lemma concerns the notion of ordinary pseudo-representation of \cite[Definition~5.3.2]{Pan22}. Let $K$ be a number field, $v$ a place of $K$ and $T:G_K\to R$ a $2$-dimensional pseudo-representation such that $T|_{G_v}=\psi_1+\psi_2$ is a sum of two characters. Recall that $T$ is \emph{$\psi_1$-ordinary} if
\[
T(\sigma\tau\eta)-\psi_1(\sigma)T(\tau\eta)-\psi_2(\tau)T(\sigma\eta)+\psi_1(\sigma)\psi_2(\tau)T(\eta)=0\qquad(\sigma,\tau\in G_v,\ \eta\in G_K).
\]
If $T$ is the trace of $\rho:G_K\to\GL_2(R)$, this means $\mathrm{tr}\big[(\rho(\sigma)-\psi_1(\sigma))(\rho(\tau)-\psi_2(\tau))\rho(\eta)\big]=0$.

\begin{lemma}\label{lem:ordline}
Let $A$ be a commutative local $\cO$-algebra with residue field $k$, and $\rho_A:G_K\to\GL_2(A)$ a representation whose reduction $\rho_k$ is absolutely irreducible. Assume that $\mathrm{tr}\rho_A|_{G_v}=\psi_1+\psi_2$ with $\psi_1(\sigma_0)-\psi_2(\sigma_0)\in A^\times$ for some $\sigma_0\in G_v$.
\begin{enumerate}
\item $\mathrm{tr}\rho_A$ is $\psi_1$-ordinary if and only if $\rho_A|_{G_v}$ has a $G_v$-stable direct summand of rank one on which $G_v$ acts by $\psi_1$. Such a summand is unique.
\item If $\mathrm{tr}\rho_A$ is both $\psi_1$-ordinary and $\psi_2$-ordinary, then $\rho_A|_{G_v}\simeq\psi_1\oplus\psi_2$.
\end{enumerate}
\end{lemma}

\begin{proof}
Since $\rho_k$ is absolutely irreducible, the $k$-span of $\rho_k(G_K)$ is $M_2(k)$ by Burnside's theorem, so the $A$-span of $\rho_A(G_K)$ is $M_2(A)$ by Nakayama's lemma. As the trace pairing on $M_2(A)$ is perfect, $\psi_1$-ordinarity is equivalent to
\[
(\rho_A(\sigma)-\psi_1(\sigma))(\rho_A(\tau)-\psi_2(\tau))=0\qquad(\sigma,\tau\in G_v).
\]
The characteristic polynomial of $\rho_A(\sigma_0)$ is $(X-\psi_1(\sigma_0))(X-\psi_2(\sigma_0))$, and its roots differ by a unit, so $A^2=L_1\oplus L_2$ with $L_i$ free of rank one and $\rho_A(\sigma_0)=\psi_i(\sigma_0)$ on $L_i$. Then $\rho_A(\sigma_0)-\psi_2(\sigma_0)$ has image $L_1$.

(1) If $\mathrm{tr}\rho_A$ is $\psi_1$-ordinary, taking $\tau=\sigma_0$ above shows that $\rho_A(\sigma)-\psi_1(\sigma)$ kills $L_1$ for all $\sigma\in G_v$, so $L_1$ is the required summand. Conversely, if $L$ is such a summand, then in a basis adapted to $L$ the product $(\rho_A(\sigma)-\psi_1(\sigma))(\rho_A(\tau)-\psi_2(\tau))$ is zero, as one checks on upper triangular matrices. Uniqueness: a stable summand with character $\psi_1$ is contained in the $\psi_1(\sigma_0)$-eigenspace of $\rho_A(\sigma_0)$, which is $L_1$.

(2) By (1), applied to $\psi_1$ and to $\psi_2$, $L_1$ and $L_2$ are both $G_v$-stable, with characters $\psi_1$ and $\psi_2$.
\end{proof}

%======================================================================
\section{Completed cohomology and ordinary parts}\label{sec:aut}

In this section $F_1$ is a totally real field, abelian over $\Q$ and of even degree $d_1$, in which $p$ splits completely. We write $\Sigma_p$ for the set of places above $p$ and fix a partition $\Sigma_p=\So\sqcup\Sn$. We use the setting of \cite[§3]{Pan22} and \cite[§§3--4]{Zha25}.

\subsection{Completed cohomology}
Let $D$ be the quaternion algebra over $F_1$ ramified exactly at the infinite places. We have $(D\otimes\mathbb A^\infty_{F_1})^\times\simeq\GL_2(\mathbb A^\infty_{F_1})$, and $D_p^\times=\prod_{v\in\Sigma_p}\GL_2(\Q_p)$. Let $U^p$ be a tame level, $\psi$ a continuous character of $(\mathbb A^\infty_{F_1})^\times/F_{1,>0}^\times$, and $\xi$ a collection of characters at the places of a finite set $\Sigma\supset\Sigma_p$, as in \cite[§4.1]{Zha25}. Completed cohomology
\[
S_{\psi,\xi}(U^p)=\varprojlim_n\varinjlim_{U_p}S_{\psi,\xi}(U^pU_p,\cO/\varpi^n)
\]
is an admissible unitary $\cO$-representation of $D_p^\times$. Let $\T=\T_{\psi,\xi}(U^p)$ be the big Hecke algebra generated by the Hecke operators $T_w$, $w\notin\Sigma$, and let $\mathfrak m$ be a maximal ideal such that $\T_{\mathfrak m}$ carries a $2$-dimensional pseudo-representation $T_{\mathfrak m}:G_{F_1,\Sigma}\to\T_{\mathfrak m}$ with $T_{\mathfrak m}(\mathrm{Frob}_w)=T_w$ and determinant $\psi\varepsilon^{-1}$, lifting $1+\bar\chi$ for a character $\bar\chi$ with $\bar\chi|_{G_v}\neq1$ for every $v\in\Sigma_p$. We write $S=S_{\psi,\xi}(U^p)_{\mathfrak m}$ and $M=\Hom_{\cO}(S\otimes_\cO E/\cO,E/\cO)$ for its Pontryagin dual, and for a prime $\mathfrak p$ of $\T_{\mathfrak m}$ we let $\rho(\mathfrak p)$ be the semisimple representation with trace $T_{\mathfrak m}\bmod\mathfrak p$.

\subsection{Paškūnas blocks}
For $v\in\Sigma_p$ let $\mathcal B_v$ be the block of the residual representation at $v$, which contains the principal series attached to $1$ and $\bar\chi|_{G_v}$. For $p=3$ and $\bar\chi|_{G_v}=\omega$ it is the exceptional block studied by Paškūnas--Tung \cite{PT21}. Let $P_{\mathcal B_v}$ be a projective generator and $E_{\mathcal B_v}=\mathrm{End}(P_{\mathcal B_v})$. Let $R^{\ps}_v$ be the pseudo-deformation ring of the residual pseudo-representation at $v$ (with fixed determinant). By \cite{Pas13} for $p\geq5$, and by \cite[Theorems~1.3 and~6.13]{PT21} in general (including the exceptional block), there is a natural map from $R^{\ps}_v$ to the centre $Z_{\mathcal B_v}$ of $E_{\mathcal B_v}$ which identifies $Z_{\mathcal B_v}$ with the maximal $\cO$-torsion-free quotient of $R^{\ps}_v$ (so $R^{\ps}_v[1/p]\simeq Z_{\mathcal B_v}[1/p]$), and $E_{\mathcal B_v}$ is finite over $R^{\ps}_v$. Only the action of $R^{\ps}_v$ through this map is used below. The local--global compatibility theorem of Pan \cite[Theorem~3.5.5]{Pan22}, extended by Zhang to the exceptional block \cite[Theorem~3.3.3]{Zha25}, states the following. Let $\mathcal B=\prod_{v}\mathcal B_v$ and $\mathbf m=\Hom(P_{\mathcal B},M)$. Then the two actions of $R^{\ps}_v$ on $\mathbf m$, through the centre of $E_{\mathcal B_v}$ and through restriction of $T_{\mathfrak m}$ to $G_v$, coincide. As a consequence, $\mathbf m$ is a faithful finitely generated $\T_{\mathfrak m}$-module \cite[Corollary~3.3.7]{Zha25}.

For a point $x$ of $R^{\ps}_v[1/p]$ we write $\mathrm{Irr}(x)$ for the set of absolutely irreducible admissible unitary Banach representations in the block which are killed by the maximal ideal of $x$.

\begin{lemma}\label{lem:irr}
Let $x$ be as above, with pseudo-representation $T_x$.
\begin{enumerate}
\item If $T_x$ is the trace of an irreducible $\rho_x$, then $\mathrm{Irr}(x)=\{\Pi_x\}$, where $\Pi_x$ corresponds to $\rho_x$ under the $p$-adic local Langlands correspondence, normalized as in \cite[\S3.3]{Zha25}, so that Colmez's functor gives $\mathbf V(\Pi_x)\simeq\rho_x\otimes\varepsilon$. If moreover $\rho_x$ is de Rham with distinct Hodge--Tate weights, then $\Pi_x$ has dense locally algebraic vectors.
\item If $T_x=\psi_1+\psi_2$ is reducible, every element of $\mathrm{Irr}(x)$ is either a character $\eta\circ\det$ or a quotient of a unitary principal series $(\Ind_{\bar B}\mu)_{\mathrm{cont}}$ for a unitary character $\mu$ of $T(\Q_p)$.
\item The absolutely irreducible admissible unitary Banach representations in $\mathcal B=\prod_v\mathcal B_v$ killed by the maximal ideal of $x=(x_v)_v$ are the completed tensor products $\widehat\bigotimes_v\Pi_v$ with $\Pi_v\in\mathrm{Irr}(x_v)$, and admissible unitary Banach representations in $\mathcal B$ killed by this maximal ideal have finite length.
\end{enumerate}
\end{lemma}

\begin{proof}
(1) The first assertion is \cite[Proposition~A.2.5]{Zha25}, that is, \cite[Proposition~6.11]{PT21}. Assume $\rho_x$ is de Rham with distinct Hodge--Tate weights. It is then potentially semistable; let $\pi$ be the smooth irreducible representation attached by the classical local Langlands correspondence to the Weil--Deligne representation of $\rho_x\otimes\varepsilon$, and $W$ the algebraic representation determined by its Hodge--Tate weights. Applying \cite[Theorem~1.3(ii)]{CDP14} to the full admissible filtered $(\varphi,N,G_{\Q_p})$-module $D_{\mathrm{pst}}(\rho_x\otimes\varepsilon)$, which lies in the set $\mathcal F(\Delta,W)$ of \cite[\S1.1]{CDP14}, we get an absolutely irreducible non-ordinary unitary completion $\Pi$ of $\pi\otimes W$ with $\mathbf V(\Pi)\simeq\rho_x\otimes\varepsilon$, and $\Pi\simeq\Pi_x$ by the injectivity in \cite[Theorem~1.1]{CDP14}. Being a completion of $\pi\otimes W$, $\Pi_x$ contains $\pi\otimes W$ as a dense subspace of locally algebraic vectors. The results of \cite{CDP14} hold for every prime $p$, in particular for $p=3$. (2) is \cite[Proposition~A.2.4]{Zha25}. (3) follows from \cite[Theorem~A.1.6 and Lemma~3.3.1]{Zha25}.
\end{proof}

\begin{lemma}\label{lem:SA}
Let $\mathfrak p$ be a maximal ideal of $\T_{\mathfrak m}[1/p]$ with $\rho(\mathfrak p)$ irreducible. Then no closed irreducible subrepresentation of $(S\otimes E)[\mathfrak p]$ has a factor of the form $\eta\circ\det$ at a place $v\in\Sigma_p$.
\end{lemma}

\begin{proof}
See \cite[proof of Proposition~4.4.6]{Pan22}: otherwise, by strong approximation for the derived group of $D^\times$, the subrepresentation would consist of functions factoring through the reduced norm, and $\rho(\mathfrak p)$ would be reducible.
\end{proof}

\subsection{Ordinary parts at the places of $\So$}
Write $G_p=D_p^\times=G_{\mathrm o}\times G_{\mathrm n}$ with $G_{\mathrm o}=\prod_{v\in\So}\GL_2(\Q_p)$ and $G_{\mathrm n}=\prod_{v\in\Sn}\GL_2(\Q_p)$. Let
\[
P=B_{\mathrm o}\times G_{\mathrm n},\qquad L_P=T_{\mathrm o}\times G_{\mathrm n},
\]
where $B_{\mathrm o}=\prod_{v\in\So}B(\Q_p)$ and $T_{\mathrm o}=\prod_{v\in\So}T(\Q_p)$. Let $N_0=\prod_{v\in\So}N(\Z_p)$. Emerton defines a functor $\Ord_P$ from smooth representations of $P$ to smooth representations of $L_P$ \cite[Definition~3.1.9]{Eme10a}, and extends it to $\varpi$-adically continuous representations \cite[Definition~3.4.1]{Eme10a}. We use the following properties.
\begin{enumerate}[label=(O\arabic*)]
\item\label{O1} $\Ord_P$ is left exact and commutes with filtered inductive limits \cite[Proposition~3.2.4]{Eme10a}.
\item\label{O2} $\Ord_P$ preserves admissibility \cite[Theorems~3.3.3 and~3.4.8]{Eme10a}.
\item\label{O3} For $V$ smooth, admissible and torsion, $\Ord_P(V)^{M_j}=e\,V^{I_{j,j}}$, the part of $V^{I_{j,j}}$ on which the Hecke operator $U_v=\sum_{a\in\Z_p/p}\left(\begin{smallmatrix}p&a\\0&1\end{smallmatrix}\right)$ acts invertibly for every $v\in\So$. Here $M_j$ is the subgroup of $T_{\mathrm o}(\Z_p)\times\prod_{v\in\Sn}\GL_2(\Z_p)$ of elements congruent to $1$ modulo $p^j$ and $I_{j,j}=\bar N_jM_jN_0$ (see \cite[proof of Theorem~3.3.3]{Eme10a}; writing $Z_M$ for the centre of $L_P$, $Z_M^+=\{z\in Z_M: zN_0z^{-1}\subset N_0\}$ and $A=\cO/\varpi^n$ with $\varpi^nV=0$, the space $V^{I_{j,j}}$ is stable under the Hecke action of $Z_M^+$ by \cite[Lemma~3.3.2]{Eme10a}, and $\Hom_{A[Z_M^+]}(A[Z_M],V^{I_{j,j}})$ is its $U_v$-invertible part by \cite[Lemma~3.1.5]{Eme10a} and its proof).
\item\label{O4} (Adjunction \cite[Theorem~4.4.6]{Eme10a}.) For $U$ admissible and $V$ smooth, or for their $\varpi$-adic versions,
\[
\Hom_{G_p}(\Ind_{\bar P}^{G_p}U,V)\simeq\Hom_{L_P}(U,\Ord_PV).
\]
By considering unit balls, the same holds for admissible unitary Banach representations.
\end{enumerate}
Put $S^{\mathrm o}=\Ord_P(S)$. For $c\geq1$ let $U(c)=U^p\cdot\prod_{v\in\Sn}K_v(c)\cdot\prod_{v\in\So}\mathrm{Iw}_1(v^c)$, where $K_v(c)$ is the principal congruence subgroup of level $p^c$ and $\mathrm{Iw}_1(v^c)$ is as in \cite[5.2.1]{Pan22}.

\medskip\noindent\textbf{Central characters.} The character $\psi$ is continuous, but in general not locally constant at $p$ (for instance $\psi=\varepsilon^2$ for the Tate module of an elliptic curve). If $\psi(z)\neq1$ for some central element $z\in U$, then $f(g)=f(gz)=\psi(z)f(g)$ forces the integral finite-level space $S_{\psi,\xi}(U,\cO)$ to vanish, although its torsion analogues need not vanish. We therefore reduce to a finite order central character, as in \cite[\S\S3.2.7, 3.3.3]{Pan22}. Let $\psi_0$ be the Teichmüller lift of $\bar\psi$. Since $\psi_0\psi^{-1}$ takes values in $1+\varpi\cO$ and $p$ is odd, it has a unique continuous square root $\theta$ with values in $1+\varpi\cO$. Multiplication by $\theta\circ\mathrm{Nrd}$ identifies $S_{\psi,\xi}(U^p)$ with $S_{\psi_0,\xi_0}(U^p)$ (after shrinking the tame level at the places of $\Sigma$, which is harmless), the action of $G_p$ being twisted by $\theta\circ\det$; modulo $\varpi^n$ this is an identification at every sufficiently deep finite level. It replaces $T_w$ by $\theta(\mathrm{Frob}_w)T_w$, hence the pseudo-representation $T$ by $\theta T$, and it commutes with $\Ord_P$, multiplying the operators $\langle\gamma\rangle_v$ defined below by $\theta_v(\Art(\gamma))$. The character $\theta$ is only used to transport identities between operators; we do not need it to be de Rham.

\begin{lemma}\label{lem:compar}
$S^{\mathrm o}/\varpi^n\simeq\Ord_P(S/\varpi^n)\simeq\varinjlim_ce\,S_{\psi,\xi}(U(c),\cO/\varpi^n)_{\mathfrak m}$, where $e=\prod_{v\in\So}\lim_mU_v^{m!}$.
\end{lemma}

\begin{proof}
The second isomorphism is \ref{O3} together with \ref{O1}. Since $S$ is $\cO$-flat and $\Ord_P$ is left exact, $\Ord_P(S)/\varpi^n$ injects into $\Ord_P(S/\varpi^n)$. For surjectivity we may replace $\psi$ by $\psi_0$, by the twist above. For $c$ large, $\psi_0$ is trivial on the centre of $U(c)$, so \cite[Corollary~3.1.3]{Pan22} applies and $S_{\psi_0,\xi_0}(U(c),\cO)\to S_{\psi_0,\xi_0}(U(c),\cO/\varpi^n)$ is surjective; so is its restriction to ordinary parts, as $e$ commutes with reduction.
\end{proof}

Let $\T^{\mathrm o}$ be the closed $\cO$-subalgebra of the continuous endomorphisms of $S^{\mathrm o}$ generated by $\T_{\mathfrak m}$ (which acts by functoriality of $\Ord_P$) and by the operators $\langle\gamma\rangle_v$, $\gamma\in\Q_p^\times$, $v\in\So$, given by the action of $\mathrm{diag}(\gamma,1)\in T(\Q_p)$ at $v$. In particular $\langle p\rangle_v=U_v$. For $v\in\So$ we define a character $\psi_{v,1}:G_v\to(\T^{\mathrm o})^\times$ by $\psi_{v,1}\circ\Art(\gamma)=\langle\gamma\rangle_v$.

\begin{lemma}\label{lem:galois}
There is a $2$-dimensional pseudo-representation $T^{\mathrm o}:G_{F_1,\Sigma}\to\T^{\mathrm o}$ with determinant $\psi\varepsilon^{-1}$ and $T^{\mathrm o}(\mathrm{Frob}_w)=T_w$ for $w\notin\Sigma$. For every $v\in\So$, we have $T^{\mathrm o}|_{G_v}=\psi_{v,1}+\psi\varepsilon^{-1}\psi_{v,1}^{-1}$, and $T^{\mathrm o}$ is $\psi_{v,1}$-ordinary.
\end{lemma}

\begin{proof}
The pseudo-representation is the push-forward of $T_{\mathfrak m}$; the content is the local identities. Since every term of these identities is multiplied by the same value of $\theta$ under the twist above, it suffices to prove them for $\psi_0$. Let $A_c$ be the $\cO$-algebra of endomorphisms of $eS_{\psi_0,\xi_0}(U(c),\cO)_{\mathfrak m}$ generated by the $T_w$ and the $\langle\gamma\rangle_v$, a finite torsion-free $\cO$-algebra. By Lemma~\ref{lem:compar}, the Hecke algebra acting on $eS_{\psi_0,\xi_0}(U(c),\cO/\varpi^n)_{\mathfrak m}$ is a quotient of $A_c/\varpi^n$, and these modules exhaust $S^{\mathrm o}/\varpi^n$; so it suffices to prove the identities in $A_c$, hence, $A_c$ being torsion-free, in $A_c[1/p]$.

We have $eS_{\psi_0,\xi_0}(U(c),E)=\bigoplus_\pi e(\pi_f^{U(c)})$, where $\pi$ runs over automorphic representations of $D^\times$ of parallel weight $2$. One-dimensional $\pi$ do not contribute, since on them $U_v$ acts by $p$ times a unit. The others are determined by their eigenvalues $T_w$, by strong multiplicity one and the Jacquet--Langlands correspondence, and for $v\in\So$ and each character of $T(\Z_p)$, the ordinary part of the corresponding isotypic component of $\pi_v^{\mathrm{Iw}_1(v^c)}$ has dimension at most one \cite[\S5.2.1]{Pan22}. Hence $A_c[1/p]$ is a product of fields. Let $\pi$ be such a representation, ordinary at $v\in\So$. By Wiles \cite[Theorem~2]{Wil88}, whose statement is place by place (see also \cite{Sai09,Ski09}), $\rho_\pi|_{G_v}$ has a subcharacter $\psi_{v,1}(\pi)$, where $\psi_{v,1}(\pi)\circ\Art$ is the character through which the $\langle\cdot\rangle_v$ act. Hence $\mathrm{tr}\rho_\pi$ has the required properties, by the first part of Lemma~\ref{lem:ordline}(1), which does not use the hypothesis on $\rho_k$.
\end{proof}

\begin{lemma}\label{lem:lgmix}
For $v\in\Sn$ and every closed $G_p$-stable subspace $\Pi\subset S\otimes E$ stable under $\T_{\mathfrak m}$, the action of $R^{\ps}_v$ on $\Ord_P(\Pi)$ through the centre of $E_{\mathcal B_v}$ coincides with its action through $T^{\mathrm o}|_{G_v}$.
\end{lemma}

\begin{proof}
An element $z$ of the centre of $E_{\mathcal B_v}$ acts on $S$ as an element $t(z)\in\T_{\mathfrak m}$, by local--global compatibility. Being a natural transformation of the identity functor on $\GL_2(\Q_p)$-representations, $z$ commutes with all $\GL_2(\Q_p)$-equivariant endomorphisms of $S$, in particular with the action of the other factors of $G_p$. Hence it is $G_p$-equivariant, and $\Ord_P(z)=t(z)$ on $S^{\mathrm o}$. On the other hand, the canonical lifting $\Ord_P(S/\varpi^n)\hookrightarrow S/\varpi^n$ \cite[Theorem~3.3.3]{Eme10a} is equivariant for the factors at $v\in\Sn$: these factors commute with $N_0$, so they lie in the monoid $M^+$ of \cite[\S3.1]{Eme10a}, and evaluation at $1$ is $M^+$-equivariant \cite[Lemma~3.1.7 and Definition~3.1.10]{Eme10a}. So the action of $z$ on $\Ord_P(S)$ through the centre is the restriction of its action on $S$, that is, $t(z)$. The same holds for $\Ord_P(\Pi)\subset S^{\mathrm o}\otimes E$, by left exactness.
\end{proof}

%======================================================================
\section{Classicality and the proof of Theorem~\ref{thm:A}}\label{sec:class}

We keep the notation of the previous section.

\subsection{Local--global compatibility at ordinary points}

\begin{proposition}\label{prop:LG}
Let $\mathfrak p$ be a maximal ideal of $\T_{\mathfrak m}[1/p]$ such that $\rho(\mathfrak p)$ is irreducible and $\rho(\mathfrak p)|_{G_v}$ is reducible for every $v\in\So$, and put $\Pi=(S\otimes E)[\mathfrak p]$. After enlarging $E$, there is a unitary character $\mu=\bigotimes_{v\in\So}\mu_v$ of $T_{\mathrm o}$, $\mu_v=\mu_{v,1}\otimes\mu_{v,2}$, such that $\Ord_P(\Pi)[T_{\mathrm o}=\mu]\neq0$ and, for every $v\in\So$, the character $\mu_{v,1}\circ\Art^{-1}$ is a subcharacter of $\rho(\mathfrak p)|_{G_v}$.
\end{proposition}

\begin{proof}
Since $\mathbf m$ is faithful and finitely generated over $\T_{\mathfrak m}$, $\Pi\neq0$, and by Lemma~\ref{lem:irr}(3) it has finite length. Let $\Pi'\subset\Pi$ be a closed irreducible subrepresentation. Then $\Pi'\simeq\widehat\bigotimes_v\Pi'_v$, and by Lemmas~\ref{lem:SA} and~\ref{lem:irr}(2), for each $v\in\So$ there are a unitary character $\mu_v$ of $T(\Q_p)$ and a surjection $f_v:(\Ind_{\bar B}\mu_v)_{\mathrm{cont}}\twoheadrightarrow\Pi'_v$. Let $U_0=\mu\otimes\widehat\bigotimes_{v\in\Sn}\Pi'_v$, a representation of $L_P$. Then $\widehat\bigotimes f_v\otimes\mathrm{id}$ is a nonzero map $\Ind_{\bar P}^{G_p}U_0\to\Pi'\subset\Pi$, and \ref{O4} gives a nonzero $L_P$-map $\phi:U_0\to\Ord_P(\Pi)$, whose image lies in $\Ord_P(\Pi)[T_{\mathrm o}=\mu]$. By \ref{O1}, $\Ord_P(\Pi)\subset S^{\mathrm o}\otimes E$. On the image of $\phi$, $\T_{\mathfrak m}$ acts through $\mathfrak p$ and $T_{\mathrm o}$ through $\mu$, so $\T^{\mathrm o}$ acts through a character, whose kernel we call $\mathfrak P$; then $\psi_{v,1}(\mathfrak P)=\mu_{v,1}\circ\Art^{-1}$. Reducing Lemma~\ref{lem:galois} modulo $\mathfrak P$, $\mathrm{tr}\rho(\mathfrak p)$ is $\psi_{v,1}(\mathfrak P)$-ordinary, and since $\rho(\mathfrak p)$ is irreducible, \cite[Lemma~5.3.3]{Pan22} (or Lemma~\ref{lem:ordline} over the residue field of $\mathfrak p$) shows that $\psi_{v,1}(\mathfrak P)$ is a subcharacter of $\rho(\mathfrak p)|_{G_v}$.
\end{proof}

\subsection{Classicality}

\begin{theorem}\label{thm:class}
Let $\mathfrak p$ be a maximal ideal of $\T_{\mathfrak m}[1/p]$, put $\Pi=(S\otimes E)[\mathfrak p]$, and let $\mu=\bigotimes_{v\in\So}\mu_v$ be a unitary character of $T_{\mathrm o}$ with $\Ord_P(\Pi)[T_{\mathrm o}=\mu]\neq0$. Assume:
\begin{enumerate}
\item $\rho(\mathfrak p)$ is irreducible;
\item for $v\in\Sn$, $\rho(\mathfrak p)|_{G_v}$ is absolutely irreducible and de Rham with distinct Hodge--Tate weights;
\item for $v\in\So$, $\rho(\mathfrak p)|_{G_v}$ is reducible and de Rham with distinct Hodge--Tate weights, and $\mu_{v,1}\circ\Art^{-1}$ is its classical character.
\end{enumerate}
Then $\rho(\mathfrak p)$ arises from a regular algebraic cuspidal automorphic representation of $\GL_2(\mathbb A_{F_1})$.
\end{theorem}

\begin{proof}
The centre of $\GL_2(\Q_p)$ at $v$ acts on $\Pi$ through the central character, which determines $\mu_{v,2}$; hence $\mu_v=\mu^{\cl}_v$ in the notation of Lemma~\ref{lem:la}. Put $U_1=\Ord_P(\Pi)[T_{\mathrm o}=\mu]\neq0$. It is an admissible unitary Banach representation of $G_{\mathrm n}$: by \ref{O2} it is admissible for $L_P$, and a compact open subgroup of $T_{\mathrm o}$ on which $\mu\equiv1\bmod\varpi$ acts trivially on the reduction of a unit ball. By Lemma~\ref{lem:lgmix}, the centre of $\mathcal B_{\mathrm n}=\prod_{v\in\Sn}\mathcal B_v$ acts on $U_1$ through the points $x_v$ defined by $\rho(\mathfrak p)|_{G_v}$, $v\in\Sn$, so by Lemma~\ref{lem:irr}(1),(3) it contains a closed subrepresentation isomorphic to $\Pi_{\mathrm n}=\widehat\bigotimes_{v\in\Sn}\Pi_{x_v}$ (any nonzero vector if $\Sn=\emptyset$). By \ref{O4} this embedding corresponds to a nonzero morphism
\[
\Phi:\widehat\bigotimes_{v\in\So}(\Ind_{\bar B}\mu^{\cl}_v)_{\mathrm{cont}}\ \widehat\otimes\ \Pi_{\mathrm n}\longrightarrow\Pi .
\]
By Lemmas~\ref{lem:la} and~\ref{lem:irr}(1), each factor of the source has dense locally algebraic vectors, hence so does the source, and $\Phi$ is nonzero on them. Hence $\Pi$ contains nonzero locally algebraic vectors, which are classical automorphic forms on $D^\times$ \cite[3.2.10]{Pan22} with Hecke eigenvalues given by $\mathfrak p$. They do not factor through the reduced norm because $\rho(\mathfrak p)$ is irreducible, and the Jacquet--Langlands correspondence gives the conclusion.
\end{proof}

\subsection{Pro-modularity}
In this subsection $\rho$ is as in Section~\ref{sec:setting} and satisfies (N), (O) and (R); we do not assume that $\rho_v$ is non-split for $v\in\So$. Let $\theta$ be the Teichmüller lift of $\bar\chi_2^{-1}$, a character of finite order, and replace $\rho$ by $\rho\otimes\theta$. This does not affect (N), (O), (R), the non-splitness of the $\rho_v$, or the modularity of $\rho$, and now $\rhobar^{ss}=1\oplus\bar\chi$. Put $\chi=\det\rho$, a de Rham character with reduction $\bar\chi$, and $\psi=\chi\varepsilon$. Since $\rho$ is odd, $\bar\chi=\det\rhobar$ is totally odd. Let $\Sigma$ be the finite set of places of $F$ consisting of the places above $p$ and those at which $\rho$ is ramified.

For a totally real field $K\supset F$ and a finite set $\Sigma_K$ of places of $K$ containing those above $\Sigma$, let $R^{\ps}_K$ be the universal deformation ring of the pseudo-representation $(1+\bar\chi)|_{G_K}$ of $G_{K,\Sigma_K}$ with determinant $\chi|_{G_K}$.

\begin{lemma}\label{lem:dimcomp}
Let $x$ be a closed point of $\Spec R^{\ps}_K[1/p]$ such that $\rho_x$ is absolutely irreducible. Then every irreducible component of $\Spec R^{\ps}_K$ containing $x$ has dimension at least $1+2[K:\Q]$.
\end{lemma}

\begin{proof}
The completed local ring of $R^{\ps}_K$ at $x$ pro-represents the functor of deformations of the pseudo-representation $T_x$ to Artinian local $k(x)$-algebras \cite[Proposition~2.2.2]{Pan22}; the proof there only uses that $G_{K,\Sigma_K}$ satisfies Mazur's finiteness condition, and $k(x)$ may be taken as coefficient field after enlarging $E$. Since $\rho_x$ is absolutely irreducible, pseudo-deformations of $T_x$ to Artinian local $k(x)$-algebras are traces of deformations of $\rho_x$, unique up to conjugation \cite{Nys96,Rou96}, so this functor coincides with the deformation functor of $\rho_x$ with fixed determinant \cite[Corollary~2.2.3]{Pan22}. Note that this uses irreducibility of $\rho_x$ over $k(x)$, not of $\bar\rho$. The latter has a presentation $k(x)[[X_1,\dots,X_{h^1}]]/(f_1,\dots,f_{h^2})$ with $h^i=\dim H^i(G_{K,\Sigma_K},\ad^0\rho_x)$. By the global Euler characteristic formula, $h^1-h^2=3[K:\Q]-[K:\Q]=2[K:\Q]$, since $H^0(G_{K,\Sigma_K},\ad^0\rho_x)=0$ and $\rho_x$ is totally odd. Hence every minimal prime of the completed local ring has coheight at least $2[K:\Q]$. Since $R^{\ps}_K$ is excellent and $x$ is a one-dimensional prime, every irreducible component through $x$ has dimension at least $2[K:\Q]+1$.
\end{proof}

\begin{lemma}\label{lem:irrF1}
Let $F'/F$ be a finite abelian extension of totally real fields. Then $\rho|_{G_{F'}}$ is absolutely irreducible.
\end{lemma}

\begin{proof}
Suppose not. Since $G_{F'}$ is normal in $G_F$, Clifford's theorem shows that $\rho|_{G_{F'}}\otimes\overline E$ is the sum of two characters $\chi_1,\chi_2$ which are conjugate under $G_F$. If $\chi_1=\chi_2$, then $\rho(G_{F'})$ consists of scalars, so the projective image of $\rho$ is finite, and for $v\mid p$ the representation $\rho_v$ is a twist of a representation with finite image; its Hodge--Tate weights are then equal, contradicting (N) or (O). If $\chi_1\neq\chi_2$, the stabilizer of $\chi_1$ is $G_M$ for a quadratic extension $M/F$ contained in $F'$, and $\rho\otimes\overline E\simeq\Ind_{G_M}^{G_F}\chi_1'$ for a character $\chi_1'$ of $G_M$ extending $\chi_1$. The character $\chi_1'$ is de Rham at the places above $p$, hence locally algebraic, hence associated with an algebraic Hecke character of $M$ \cite[Chapter~III]{Ser68}. Since $M$ is totally real, $\chi_1'$ is the product of a power $\varepsilon^n$ of the cyclotomic character and a character of finite order. Both Jordan--Hölder constituents of $\rho_v|_{G_{M_w}}$, for $w\mid v\mid p$, are then of the form $\varepsilon^n\cdot(\text{finite order})$, so the Hodge--Tate weights of $\rho_v$ are equal, a contradiction.
\end{proof}

\begin{proposition}\label{prop:promod}
There is a totally real field $F_1\supset F$, abelian over $\Q$, of even degree and with $p$ totally split in $F_1$, with the following property. Let $\Sigma_1$ be the set of places of $F_1$ above $\Sigma$. Every irreducible component $C$ of $\Spec R^{\ps}_{F_1}$ (with $\Sigma_{F_1}=\Sigma_1$) containing the point defined by $\rho|_{G_{F_1}}$ is pro-modular in the sense of \cite[Definition~4.1.3]{Zha25}: there are a tame level $U^p$, a finite set of places $\Sigma'\supset\Sigma_1$ and characters $\xi$ as in \cite[§4.1]{Zha25} such that $C$ lies in the image of the closed immersion $\Spec\T_{\psi,\xi}(U^p)_{\mathfrak m}\to\Spec R^{\ps}_{F_1,\Sigma'}$, where $R^{\ps}_{F_1,\Sigma'}$ is the pseudo-deformation ring of $G_{F_1,\Sigma'}$, of which $R^{\ps}_{F_1}$ is a quotient.
\end{proposition}

\begin{proof}
Take $F_1$ as in \cite[Theorem~5.5.1]{Zha25}, applied with the set $\Sigma$. Its hypotheses are satisfied by (R): $\bar\chi$ is totally odd, extends to $G_\Q$, and $\bar\chi|_{G_v}\neq1$ for $v\mid p$. By Lemma~\ref{lem:irrF1}, $\rho|_{G_{F_1}}$ is absolutely irreducible, so by Lemma~\ref{lem:dimcomp}, $C$ has dimension at least $1+2[F_1:\Q]$. Hence it is pro-modular by \cite[Theorem~5.5.1]{Zha25}, and by \cite[Remark~4.1.2]{Zha25} a single pair $(U^p,\xi)$ works for the whole of $C$. The proof of this theorem does not use any local condition at $p$ on the points of $C$.
\end{proof}

From now on we also write $\So$ (resp.\ $\Sn$) for the set of places of $F_1$ above the places of $\So$ (resp.\ $\Sn$). Since $p$ splits completely in $F_1$, for a place $v$ of $F_1$ above $u\in\Sigma_p$ we have $G_v=G_{F_u}$ and $\rho|_{G_v}=\rho_u$.

\subsection{Proof of Theorem~\ref{thm:A}}
Let $F_1$ be as in Proposition~\ref{prop:promod} and $C$ a component through the point of $\rho|_{G_{F_1}}$. With $(U^p,\xi)$ as there, $\rho|_{G_{F_1}}$ defines a maximal ideal $\mathfrak p$ of $\T_{\mathfrak m}[1/p]$, and $\rho(\mathfrak p)=\rho|_{G_{F_1}}$ is irreducible by Lemma~\ref{lem:irrF1}. By Proposition~\ref{prop:LG}, there is $\mu$ with $\Ord_P((S\otimes E)[\mathfrak p])[T_{\mathrm o}=\mu]\neq0$ such that $\mu_{v,1}\circ\Art^{-1}$ is a subcharacter of $\rho|_{G_v}$ for $v\in\So$. As $\rho|_{G_v}$ is non-split, this is its unique subcharacter, which is the classical one by Lemma~\ref{lem:BK}. Theorem~\ref{thm:class} shows that $\rho|_{G_{F_1}}$ is modular. Since $F_1/F$ is abelian and $\rho|_{G_{F_1}}$ is irreducible, $\rho$ is modular by solvable descent \cite[Lemma~2.2.2]{BLGGT14}, \cite{AC89}. \qed

%======================================================================
\section{Three applications of Theorem~\ref{thm:A}}\label{sec:appA}
We say that a representation $r:G_{\Q_p}\to\GL_2(E)$ is \emph{of Steinberg type} if $r\otimes\eta$ is semistable and not crystalline for some de Rham character $\eta$ of $G_{\Q_p}$; equivalently, $r$ is de Rham and its Weil--Deligne representation is a twist of the one attached to the Steinberg representation. We say that $r$ is \emph{of supercuspidal type} if it is de Rham and its Weil--Deligne representation (attached to $D_{\mathrm{pst}}(r)$) is irreducible, and \emph{of discrete series type} if it is of Steinberg or of supercuspidal type, that is, if its Frobenius-semisimple Weil--Deligne representation $(r_{\mathrm{WD}},N)$ is indecomposable (either irreducible, or reducible with $N\neq0$).

\begin{lemma}\label{lem:St}
Let $r:G_{\Q_p}\to\GL_2(E)$.
\begin{enumerate}
\item If $r$ is of Steinberg type, then $r$ has distinct Hodge--Tate numbers. If they differ by $1$, then $r$ is reducible and non-split; otherwise $r$ is irreducible.
\item If $r$ is of supercuspidal type, then $r$ is irreducible.
\end{enumerate}
\end{lemma}

\begin{proof}
(1) Twisting, we may assume that $r$ is semistable and not crystalline. Let $D=D_{\mathrm{st}}(r)$, a weakly admissible filtered $(\varphi,N)$-module of rank $2$ over $E\otimes\Q_p=E$ \cite{CF00}, with $N\neq0$. Choose a basis $e_1,e_2$ with $Ne_2=e_1$ and $Ne_1=0$. Since $N\varphi=p\varphi N$, the line $Ee_1=\ker N$ is $\varphi$-stable, say $\varphi(e_1)=\alpha e_1$, and $\varphi(e_2)\in p\alpha e_2+Ee_1$. Hence $t_N(D)=2v_p(\alpha)+1$. Let $h_1\leq h_2$ be the jumps of the Hodge filtration, so that $t_H(D)=h_1+h_2$ and weak admissibility gives $v_p(\alpha)=(h_1+h_2-1)/2$. If $h_1=h_2$ the filtration is trivial and $t_H(Ee_1)=h_1>v_p(\alpha)$, which contradicts weak admissibility; hence $h_1<h_2$. A $(\varphi,N)$-stable line is contained in $\ker N$, so $Ee_1$ is the only one, and $r$ is reducible if and only if $Ee_1$ with the induced filtration is weakly admissible, that is, $t_H(Ee_1)=v_p(\alpha)$. As $t_H(Ee_1)\in\{h_1,h_2\}$ and $v_p(\alpha)=(h_1+h_2-1)/2$, this forces $h_2=h_1+1$; conversely, if $h_2=h_1+1$ then $v_p(\alpha)=h_1$, and weak admissibility of $D$ gives $t_H(Ee_1)\leq h_1$, so $t_H(Ee_1)=h_1=v_p(\alpha)$. In that case $r$ is non-split, since a sum of two semistable characters has $N=0$.

(2) If $r$ were reducible, an extension of $\psi_2$ by $\psi_1$, then by the exactness of $D_{\mathrm{pst}}$ its Weil--Deligne representation would contain that of $\psi_1$, contradicting irreducibility.
\end{proof}

\begin{remark}\label{rem:St}
The argument in (2) is the same as in \cite[Theorem~3.2]{JK}, where it is used to show that the restriction to a decomposition group at $p$ of the $p$-adic representation attached to a newform $f$ is irreducible when $\pi_p(f)$ is supercuspidal. Part (1) is well known: a two-dimensional semistable non-crystalline representation of $G_{\Q_p}$ is reducible (and then non-split) if and only if its Hodge--Tate numbers differ by $1$, see \cite{Bre03}, \cite{BM02}.
\end{remark}

\begin{corollary}\label{cor:St}
Let $p$, $F$ and $\rho$ be as in Section~\ref{sec:setting}, and assume that $\rho$ satisfies \textup{(R)} and that for every $v\in\Sigma_p$, $\rho_v$ has distinct Hodge--Tate numbers and is of discrete series type. Then $\rho$ is modular.
\end{corollary}

\begin{proof}
Let $\So$ be the set of $v\in\Sigma_p$ at which $\rho_v$ is of Steinberg type with Hodge--Tate numbers differing by $1$, and $\Sn=\Sigma_p\setminus\So$. By Lemma~\ref{lem:St}, $\rho_v$ is irreducible for $v\in\Sn$, and reducible and non-split for $v\in\So$, so (N) and (O) hold. Theorem~\ref{thm:A} applies.
\end{proof}

For the remaining applications, let $F$ be a real quadratic field and $p$ an odd prime which splits in $F$, with $\Sigma_p=\{v_1,v_2\}$.

\begin{corollary}\label{cor:St2}
Let $\rho:G_F\to\GL_2(E)$ be a continuous, absolutely irreducible, odd representation, unramified outside finitely many places, satisfying \textup{(R)}. Assume that for $v\in\Sigma_p$, $\rho_v$ is de Rham with Hodge--Tate numbers differing by $1$, and that $\rho_{v_s}$ is of Steinberg type for some $v_s\in\Sigma_p$. Then $\rho$ is modular.
\end{corollary}

\begin{proof}
By Lemma~\ref{lem:St}, $\rho_{v_s}$ is reducible and non-split. Let $v$ be the other place above $p$. If $\rho_v$ is irreducible, Theorem~\ref{thm:A} applies with $\So=\{v_s\}$. If $\rho_v$ is reducible, then $\rho$ is reducible at both places above $p$, and the classical character of $\rho_w$ (Lemma~\ref{lem:BK}) is a subcharacter with strictly smaller Hodge--Tate number for $w\in\Sigma_p$; by (R), after twisting by a character of finite order, $\rho$ satisfies the hypotheses of \cite[Theorem~5.1.2]{Pan22}, which allows $\rho_v$ to be split, and $\rho$ is modular.
\end{proof}

\begin{theorem}\label{thm:2k}
Let $\rho:G_F\to\GL_2(E)$ be a continuous, absolutely irreducible, odd representation, unramified outside finitely many places, satisfying \textup{(R)}. Assume:
\begin{enumerate}
\item for $v\in\Sigma_p$, $\rho_v$ is de Rham with distinct Hodge--Tate numbers; these differ by $1$ at one place $v_c$ and by $k-1\geq2$ at the other place $v_k$;
\item $\rho$ is potentially modular: there are a finite totally real Galois extension $W/F$, a regular algebraic cuspidal automorphic representation $\pi$ of $\GL_2(\mathbb A_W)$ and an isomorphism $\iota:\overline\Q_p\simeq\mathbb C$ such that $\rho|_{G_W}\simeq r_{\pi,\iota}$;
\item $\rho_{v_s}$ is of Steinberg type for some $v_s\in\Sigma_p$.
\end{enumerate}
Then:
\begin{enumerate}[label=\textup{(\alph*)}]
\item There is a continuous representation $\rho':G_F\to\GL_2(E')$, with $E'/\Q_p$ finite, which is absolutely irreducible, odd and unramified outside finitely many places, with the following properties. After extending scalars, $\rho$ and $\rho'$ are the members at two places $\lambda_0,\lambda'$ above $p$ of a strongly compatible system with coefficients in a number field $M$; thus for every finite place $w$ of $F$ their Frobenius-semisimple Weil--Deligne representations (via $D_{\mathrm{pst}}$ if $w\mid p$) are obtained from a common Weil--Deligne representation over $\overline M$ through the embeddings $\overline M\to\overline\Q_p$ defined by $\lambda_0$ and $\lambda'$. In particular $\rho_w$ is of Steinberg type if and only if $\rho'_w$ is. The representation $\rho'_v$ is de Rham for $v\in\Sigma_p$, the Hodge--Tate numbers of $\rho'_{v_1}$ (resp.\ $\rho'_{v_2}$) being those of $\rho_{v_2}$ (resp.\ $\rho_{v_1}$), and $\rho$ is modular if $\rho'$ is.
\item If $v_s=v_c$, then $\rho$ is modular.
\item If $v_s=v_k$, assume that either $\rhobar'$ is reducible and $\rho'$ satisfies \textup{(R)}, or $\rhobar'$ is absolutely irreducible and modular and, if $p=3$, $\rhobar'|_{G_{F(\zeta_3)}}$ is absolutely irreducible. Then $\rho$ is modular.
\end{enumerate}
\end{theorem}

Note that $k$ is necessarily even: $\det\rho$ is a de Rham character of the totally real field $F$, so it is a power of $\varepsilon$ times a character of finite order \cite[Chapter~III]{Ser68}, and the sums of the two Hodge--Tate numbers at $v_c$ and at $v_k$ coincide; hence $1$ and $k-1$ have the same parity.

Strictly speaking, the isomorphisms in (a) are not isomorphisms of $\overline\Q_p$-representations: if the Weil--Deligne representations of $\rho$ and $\rho'$ were isomorphic as such at every place, $\rho$ and $\rho'$ would have the same Frobenius traces, hence would be isomorphic, which is impossible since their Hodge--Tate numbers differ.

\begin{proof}
(a) We use the descent construction of \cite{Die04} and of the proof of \cite[Theorem~5.5.1]{BLGGT14}, starting from the potential modularity given by (2); we do not use the residual hypotheses in the statement of \cite[Theorem~5.5.1]{BLGGT14}. Solvable descent gives, for the intermediate fields $W\supset W_i\supset F$ with $W/W_i$ solvable, cuspidal representations whose $p$-adic representations are $\rho|_{G_{W_i}}$, and Brauer's induction theorem combines their compatible systems into a strongly compatible system $(\rho_\lambda)_\lambda$ of representations of $G_F$ with coefficients in a number field $M$, indexed by the finite places $\lambda$ of $M$, with $\rho\otimes\overline\Q_p=\rho_{\lambda_0}\otimes\overline\Q_p$ for some $\lambda_0\mid p$. Compatibility at the places above the residue characteristic follows from local--global compatibility at $p$ for Hilbert modular forms \cite{Sai09}, \cite[Theorem~1]{Ski09}. As in \cite[\S5.1]{BLGGT14}, the Hodge--Tate numbers of the members of the system are given by multisets $H_\tau$ indexed by the embeddings $\tau:F\to\overline M$, which are invariant under $\Gal(\overline M/M)$: for $\lambda\mid\ell$ and $v\mid\ell$, the Hodge--Tate numbers of $\rho_\lambda|_{G_{F_v}}$ are $H_\tau$, where $\tau$ is the embedding $F\to F_v\to\overline M_\lambda$ restricted to an embedding into $\overline M$. By (1), $H_{\tau_1}\neq H_{\tau_2}$ for the two embeddings $\tau_1,\tau_2$ of $F$, so they are not conjugate under $\Gal(\overline M/M)$, and $M$ contains $\tau_1(F)$. Identify $F$ with $\tau_1(F)\subset M$. For a place $\lambda\mid p$ of $M$, let $u(\lambda)$ be the place of $F$ below $\lambda$; then $\rho_\lambda|_{G_{F_{u(\lambda)}}}$ has Hodge--Tate numbers $H_{\tau_1}$, and $\rho_\lambda$ has Hodge--Tate numbers $H_{\tau_2}$ at the other place above $p$. Since $p$ splits in $F$, there is a place $\lambda'\mid p$ of $M$ with $u(\lambda')\neq u(\lambda_0)$, and we put $\rho'=\rho_{\lambda'}$. It is unramified outside finitely many places and de Rham at the places above $p$, its Weil--Deligne representations are related to those of $\rho$ as stated, by strong compatibility, and its Hodge--Tate numbers at $v_1$ and $v_2$ are interchanged with respect to those of $\rho$. Moreover $\rho'|_{G_W}\simeq r_{\pi,\iota'}$ for a suitable $\iota':\overline\Q_p\simeq\mathbb C$, so $\rho'$ is odd, as $W$ is totally real, and absolutely irreducible. Finally, suppose that $\rho'\otimes\overline\Q_p\simeq r_{\pi',\iota''}$ for a regular algebraic cuspidal $\pi'$ of $\GL_2(\mathbb A_F)$ and $\iota'':\overline\Q_p\simeq\mathbb C$. Choose $\iota:\overline\Q_p\simeq\mathbb C$ such that $\iota\circ j_0=\iota''\circ j'$ on $M$, where $j_0,j':M\to\overline\Q_p$ are the embeddings defined by $\lambda_0$ and $\lambda'$. Then $\rho$ and $r_{\pi',\iota}$ have the same Frobenius traces, so $\rho\otimes\overline\Q_p\simeq r_{\pi',\iota}$, as $\rho$ is irreducible.

(b) The proof of Corollary~\ref{cor:St2} applies verbatim with $v_s=v_c$: it only uses that the Hodge--Tate numbers of $\rho_{v_s}$ differ by $1$ and that those of the other place are distinct.

(c) By Lemma~\ref{lem:St}, $\rho_{v_k}$ is irreducible. Let $\rho'$ be as in (a); it suffices to show that $\rho'$ is modular. By (a), $\rho'_{v_k}$ is of Steinberg type, and its Hodge--Tate numbers, those of $\rho_{v_c}$, differ by $1$; by Lemma~\ref{lem:St}, $\rho'_{v_k}$ is reducible and non-split. The Hodge--Tate numbers of $\rho'_{v_c}$ are distinct. If $\rhobar'$ is reducible and $\rho'$ satisfies (R), then, as in (b), $\rho'$ is modular by Theorem~\ref{thm:A} with $\So=\{v_k\}$ if $\rho'_{v_c}$ is irreducible, and by \cite[Theorem~5.1.2]{Pan22} if $\rho'_{v_c}$ is reducible. Suppose now that $\rhobar'$ is absolutely irreducible and modular. We claim that $\rhobar'|_{G_{F(\zeta_p)}}$ is absolutely irreducible; for $p=3$ this is part of the hypothesis, so let $p>3$. Since $p$ is unramified in $F$, $\Gal(F(\zeta_p)/F)\simeq(\Z/p\Z)^\times$ is cyclic, and $K=F(\sqrt{p^*})$, with $p^*=(-1)^{(p-1)/2}p$, is its unique quadratic subextension. The proof of \cite[Lemma~1.13]{DP23}, written over $\Q$, applies verbatim over $F$ and shows that if $\rhobar'|_{G_{F(\zeta_p)}}$ is not absolutely irreducible, then $\rhobar'$ is bad dihedral, that is, induced from a character of $G_K$. In that case, since $K/F$ is ramified at $v_k$ and $F_{v_k}=\Q_p$, the restriction $\rhobar'|_{G_{F_{v_k}}}$ is induced from a character of $G_{\Q_p(\sqrt{p^*})}$. The proof of \cite[Lemma~1.14]{DP23}, which follows \cite[Proposition~2.2]{Rib97}, only uses this local information at $p$, and therefore applies to $\rhobar'|_{G_{F_{v_k}}}$: if $k_0$ is its Serre weight, then $p=2k_0-1$ or $p=2k_0-3$. Since this argument only involves the projective image, we may twist $\rho'$ by a character, and since $\rho'_{v_k}$ is of Steinberg type with Hodge--Tate numbers differing by $1$, we may assume that $\rho'_{v_k}$ is semistable with Hodge--Tate numbers $0$ and $1$, so that $k_0\in\{2,p+1\}$. The projective image of a bad dihedral representation is dihedral of order prime to $p$, so the wild inertia subgroup at $v_k$ acts trivially and $k_0\neq p+1$. Hence $k_0=2$ and $p=3$, a contradiction. This proves the claim, and then $\rho'$ is modular by the results of Kisin \cite{Kis09}, which hold over totally real fields in which $p$ splits completely and allow arbitrary de Rham local conditions with distinct Hodge--Tate weights at the places above $p$, as extended by Hu--Tan \cite{HT15} and Tung \cite{Tun21} to all odd primes.
\end{proof}

\begin{remark}\label{rem:entangled}
Theorem~\ref{thm:2k} may look somewhat technical, but it applies to representations of weight $(2,k)$ with ramification of Steinberg type at one of the two places above $p$ (the case of parallel weight $2$ being covered by Corollary~\ref{cor:St2}), and we would like to stress its conceptual content. Using the hypothesis of potential modularity, a second representation $\rho'$ is constructed from $\rho$ by a suitable conjugation inside the compatible system containing $\rho$: its Weil--Deligne representations at all finite places are those of $\rho$, up to the change of coefficient embedding, and its Hodge--Tate numbers at the two places above $p$ are interchanged with respect to those of $\rho$. We call $(\rho,\rho')$ an \emph{entangled pair}. The two members of the pair are members of the same compatible system, the only difference being the place of the coefficient field above $p$ at which the $p$-adic representation is taken, so their residual representations are in general unrelated, but the modularity of one of them implies that of the other (see (a) and its proof). The hypotheses of Theorem~\ref{thm:2k}\textup{(c)} are then the standard hypotheses of modularity lifting theorems of Wiles type, imposed on the two members of the pair: residual reducibility together with condition (R), or residual absolute irreducibility together with residual modularity (the usual condition of large residual image being automatic here for $p>3$, by the local argument in the proof). In this way the modularity of $\rho$, a residually reducible mixed case (if $\rho$ is not mixed, the results of Pan and Zhang already apply), is reduced either to a residually reducible case covered by Theorem~\ref{thm:A} and \cite[Theorem~5.1.2]{Pan22}, or to a residually irreducible case, where modularity lifting theorems covering all cases, including the mixed ones, are available \cite{Kis09,HT15,Tun21}. Theorem~\ref{thm:2k} is therefore a modularity lifting theorem of Wiles type for a pair of entangled representations. Finally, this type of result is meant to be applied in inductive arguments, where modularity is proved by induction on the size of invariants such as the weight and the level, and is propagated through safe chains of congruences. In such arguments potential modularity is available, and the entangled representation $\rho'$ has invariants as small as those of $\rho$: the same level, since its Weil--Deligne representations at all finite places are those of $\rho$ up to the change of coefficient embedding, and the same Hodge--Tate numbers, interchanged between the two places above $p$. Therefore, when $\rhobar'$ is irreducible, its modularity is provided by the induction hypothesis, as is inherent to the inductive argument, and the hypotheses of \textup{(c)} are not restrictive in these applications (see Remark~\ref{rem:applications}).
\end{remark}

%======================================================================
\section{Abelian varieties of $\GL_2$-type}\label{sec:AV}

In this section $p$ and $F$ are as in Section~\ref{sec:setting}. Let $A/F$ be an abelian variety of dimension $g$ and $K\subset\mathrm{End}_F(A)\otimes\Q$ a number field with $[K:\Q]=g$. Replacing $A$ by an isogenous abelian variety over $F$, we may and do assume that $\cO_K\subset\mathrm{End}_F(A)$. For a prime $\lambda$ of $K$ above $p$ let $T_\lambda A$ be the $\lambda$-adic Tate module and $V_\lambda A=T_\lambda A\otimes_{\cO_{K,\lambda}}K_\lambda$.

\begin{lemma}\label{lem:AVbasic}
\begin{enumerate}
\item $V_\lambda A$ is free of rank $2$ over $K_\lambda$; we write $\rho_{A,\lambda}:G_F\to\GL_2(K_\lambda)$ for the corresponding representation.
\item $\rho_{A,\lambda}$ is unramified outside the places above $p$ and those of bad reduction, it is semisimple, and for $v\mid p$ it is de Rham with Hodge--Tate numbers $0$ and $-1$.
\item $\rho_{A,\lambda}$ is totally odd.
\end{enumerate}
\end{lemma}

\begin{proof}
(1) is due to Ribet \cite{Rib76b,Rib92}. In (2), unramifiedness is the criterion of Néron--Ogg--Shafarevich \cite{ST68}, semisimplicity is due to Faltings \cite{Fal83}, and the de Rham property, with the Hodge--Tate decomposition of $V_pA$, is due to Faltings \cite{Fal88}. For the Hodge--Tate numbers of each coefficient factor, decompose $H_1(A_\sigma(\mathbb C),\mathbb C)=\bigoplus_\tau H_\tau$ along the embeddings $\tau:K\to\mathbb C$; each $H_\tau$ has dimension $2$, and complex conjugation, which commutes with $K$, exchanges $H^{-1,0}$ and $H^{0,-1}$, so each $H_\tau$ meets both in dimension one. (3) Let $c$ be the complex conjugation attached to a real place $\sigma$ of $F$. It acts on $H_1(A_\sigma(\mathbb C),\Q)$, a $K$-vector space of dimension $2$, as a $K$-linear involution whose $\pm1$-eigenspaces have $\Q$-dimension $g$ each, since $c$ exchanges $H^{-1,0}$ and $H^{0,-1}$. Hence each eigenspace has $K$-dimension $1$, and $\det\rho_{A,\lambda}(c)=-1$ by the comparison isomorphism $H_1(A_\sigma(\mathbb C),\Q)\otimes_KK_\lambda\simeq V_\lambda A$.
\end{proof}

\begin{lemma}\label{lem:AVirr}
$\rho_{A,\lambda}$ is absolutely irreducible.
\end{lemma}

\begin{proof}
This is well known (see \cite{Rib76b}); we recall the argument. Suppose that $\rho_{A,\lambda}\otimes\overline{K}_\lambda$ has a $G_F$-stable line, with character $\eta$. Then $\eta$ is de Rham at the places above $p$, hence locally algebraic and associated with an algebraic Hecke character of $F$ \cite[Chapter~III]{Ser68}. Since $F$ is totally real, $\eta=\varepsilon^n\eta_0$ with $n\in\Z$ and $\eta_0$ of finite order. For a place $w\nmid p$ of good reduction, $\eta(\mathrm{Frob}_w)$ is then, up to a root of unity, $N(w)^{\pm n}$, while by Weil the eigenvalues of $\mathrm{Frob}_w$ on $V_\lambda A$ have absolute value $N(w)^{1/2}$ under every complex embedding. This is impossible.
\end{proof}

\subsection{The CM case}

\begin{lemma}\label{lem:AVCM}
Assume that $A_{\overline\Q}$ has complex multiplication. Then $\rho_{A,\lambda}$ is modular.
\end{lemma}

\begin{proof}
This is classical. By the theory of complex multiplication \cite{ST61}, $\rho_{A,\lambda}$ becomes abelian over a finite extension of $F$, so the identity component of the Zariski closure of its image is a torus, non-central since the Hodge--Tate numbers are distinct. As $\rho_{A,\lambda}$ is absolutely irreducible (Lemma~\ref{lem:AVirr}), $G_F$ permutes the two weight lines of this torus nontrivially, so $\rho_{A,\lambda}\otimes\overline K_\lambda\simeq\Ind_{G_M}^{G_F}\theta$ for a quadratic extension $M/F$ and a character $\theta\neq\theta^c$ of $G_M$. The character $\theta$ is associated with an algebraic Hecke character of $M$ \cite[Chapter~III]{Ser68}, which is pure of weight $1$ by Weil's bounds for $A$. At a real place of $M$ its infinity type would be an integer $a$ with $2a=1$; hence $M$ has no real place and is a CM field. The theta series of $\theta$ (equivalently, its automorphic induction) is cuspidal since $\theta\neq\theta^c$, and gives a Hilbert eigenform of parallel weight $2$ over $F$ whose Galois representation is $\rho_{A,\lambda}$ \cite{Hec20,AC89}.
\end{proof}

\subsection{Local indecomposability}

\begin{lemma}\label{lem:AVnonsplit}
Assume that $A_{\overline\Q}$ does not have complex multiplication, and let $v\mid p$ be such that $\rho_v=\rho_{A,\lambda}|_{G_{F_v}}$ is reducible. Then $\rho_v$ is non-split.
\end{lemma}

\begin{proof}
It suffices to show that $\rho_{A,\lambda}|_{G_{k_w}}$ is non-split for some finite extension $k_w$ of $F_v$. By Grothendieck's semistable reduction theorem \cite[Exp.~IX]{SGA7}, there is a finite extension $k/F$ of number fields and a place $w\mid v$ of $k$ such that $A_k$ has semistable reduction at $w$; we may take $k=F(A[N])$ for some $N\geq3$ prime to $p$. Let $t$ be the toric rank of the special fibre of the Néron model of $A_k$ at $w$. The character group of the torus is a free $\Z$-module of rank $t$ on which $\cO_K$ acts, so $t$ is divisible by $[K:\Q]=g$, and $t\in\{0,g\}$.

If $t=g$, the reduction is purely toric. By Grothendieck's monodromy theorem \cite[Exp.~IX]{SGA7}, $\rho_{A,\lambda}|_{G_{k_w}}$ is semistable with non-zero monodromy operator. A split representation would be a sum of two semistable characters, which are crystalline since the monodromy of a character vanishes. Hence $\rho_{A,\lambda}|_{G_{k_w}}$ is non-split.

If $t=0$, $A_k$ has good reduction at $w$, and $\rho_{A,\lambda}|_{G_{k_w}}$ is crystalline and reducible. The Jordan--Hölder characters of $\rho_v$ have Hodge--Tate numbers $0$ and $-1$, and a stable line, a priori defined over a finite extension of $K_\lambda$, descends to $K_\lambda$, since its Galois conjugates are stable lines with the same Hodge--Tate number. These characters are characters of $G_{F_v}=G_{\Q_p}$, so after a finite extension of $k_w$ they are an unramified character and $\varepsilon$ times an unramified character. Hence the factor $V_\lambda A$ of $V_pA$ is ordinary: the corresponding $\lambda$-divisible group has slopes $0$ and $1$.

We now pass to real multiplication, keeping track of this factor. Since $A_{\overline\Q}$ has no complex multiplication and $A$ is defined over the totally real field $F$, by \cite[Propositions~1.1 and~1.4]{Zhao14} $A_{\overline\Q}$ is isogenous to a power $B^e$ of an absolutely simple abelian variety $B$, and $\mathrm{End}^0(B)$ contains a totally real field $L$ of degree $\dim B$. Enlarging $k$, we may assume that $B$, $L$ and the isogeny are defined over $k$, and after a further isogeny that $\cO_L$ acts on $B$. As in the proof of \cite[Theorem~5.6]{Zhao14}, comparing the actions of $K$ and $L$ through the centre of the simple algebra $\mathrm{End}^0(A_k)$, we obtain, after a finite extension $E'$ of the coefficients, an isomorphism $V_\lambda A\otimes_{K_\lambda}E'\simeq V_\mu B\otimes_{L_\mu}E'$ of representations of $G_{k_w}$ for some $\mu\mid p$. In particular the $\mu$-factor of $B$ at $w$ is ordinary, and \cite[Theorem~5.2]{Zhao14} shows that $V_\mu B$ is indecomposable as a representation of the inertia group at $w$. A nonzero extension class remains nonzero after extending coefficients, so $\rho_{A,\lambda}|_{G_{k_w}}$, and hence $\rho_v$, is non-split.
\end{proof}

\begin{remark}
For elliptic curves Lemma~\ref{lem:AVnonsplit} is classical: if $E$ has good ordinary reduction and $\rho_{E,p}|_{G_{k_w}}$ is split, then its Serre--Tate parameter is a root of unity, so $E$ is a quasi-canonical lift of its reduction (isogenous to the canonical lift), and it has complex multiplication. For abelian varieties with real multiplication it is due to Hida \cite{Hid13} when $p$ is unramified in $k$ and in $K$, and to Zhao \cite{Zhao14} in general. Zhao's theorem is a statement about abelian varieties and does not assume modularity.
\end{remark}

\subsection{Proof of Theorem~\ref{thm:B}}
Write $\rho=\rho_{A,\lambda}$, and let $E=K_\lambda$. If $A_{\overline\Q}$ has complex multiplication, $\rho$ is modular by Lemma~\ref{lem:AVCM}. Assume it does not. By Lemmas~\ref{lem:AVbasic} and~\ref{lem:AVirr}, $\rho$ is continuous, absolutely irreducible, odd and unramified outside finitely many places, and for all $v\mid p$ it is de Rham with distinct Hodge--Tate weights. After extending scalars as in Section~\ref{sec:setting}, let $\So$ be the set of $v\mid p$ such that $\rho_v$ is reducible, and $\Sn$ its complement. Hence (N) and (O) hold, (R) holds by assumption, and $\rho_v$ is non-split for $v\in\So$ by Lemma~\ref{lem:AVnonsplit}. By Theorem~\ref{thm:A}, $\rho$ is modular. Since its Hodge--Tate numbers are $0$ and $-1$ at every $v\mid p$, the corresponding automorphic representation has parallel weight $2$. \qed

\subsection{Potentially modular representations of parallel weight $2$}\label{sec:Bprime}

\begin{thmBprime}
Let $\rho:G_F\to\GL_2(E)$ be continuous, absolutely irreducible and unramified outside finitely many places, and assume:
\begin{enumerate}
\item for every $v\mid p$, $\rho_v$ is de Rham with Hodge--Tate numbers $0$ and $-1$;
\item $\rho$ satisfies \textup{(R)};
\item $\rho$ is potentially modular: there are a finite totally real Galois extension $W/F$, a regular algebraic cuspidal automorphic representation $\pi$ of $\GL_2(\mathbb A_W)$ and an isomorphism $\iota:\overline\Q_p\simeq\mathbb C$ such that $\rho|_{G_W}\simeq r_{\pi,\iota}$;
\item either $[F:\Q]$ is odd, or there is a finite place $v_0$ of $F$ such that the Weil--Deligne representation attached to $\rho|_{G_{F_{v_0}}}$ is indecomposable.
\end{enumerate}
Then $\rho$ is modular.
\end{thmBprime}

In (4), if $v_0\mid p$ the Weil--Deligne representation is the one attached to $D_{\mathrm{pst}}(\rho_{v_0})$. A two-dimensional Frobenius-semisimple Weil--Deligne representation is indecomposable if and only if it corresponds under the local Langlands correspondence to a discrete series representation, that is, to a supercuspidal representation or to a twist of the Steinberg representation. For instance, (4) holds if $\rho|_{G_{F_v}}$ has supercuspidal type at some finite place $v$ of $F$.

\begin{proof}
By Theorem~\ref{thm:B}, it suffices to show that $\rho\simeq\rho_{A,\lambda}$ for an abelian variety $A/F$ of $\GL_2$-type and a prime $\lambda\mid p$ of its field of coefficients, after extending scalars. Since $W$ is totally real and $\rho|_{G_W}\simeq r_{\pi,\iota}$, the representation $\rho$ is odd, and by (1) $\pi$ has parallel weight $2$.

Assume first that $\pi$ has complex multiplication, so that $\rho|_{G_W}$ is induced from a character of a CM quadratic extension $M_W$ of $W$. Then the identity component of the Zariski closure of $\rho(G_F)$ is a torus, non-central since the Hodge--Tate numbers are distinct, and as $\rho$ is absolutely irreducible, $G_F$ permutes its two weight lines nontrivially. Thus $\rho\otimes\overline E\simeq\Ind_{G_M}^{G_F}\theta$ for a quadratic extension $M/F$ and a character $\theta\neq\theta^c$ of $G_M$, and the stabilizer of the lines in $G_W$ gives $MW=M_W$. If $M$ had a real embedding, it would extend to a real embedding of $MW$, as $W$ is totally real; so $M$ is totally imaginary, hence CM. The character $\theta$ is associated with an algebraic Hecke character of $M$ \cite[Chapter~III]{Ser68}, and $\rho$ is modular as in the proof of Lemma~\ref{lem:AVCM}.

Assume now that $\pi$ does not have complex multiplication. The descent construction recalled in the proof of Theorem~\ref{thm:2k}\textup{(a)}, which does not require $F$ to be quadratic (see also \cite[\S9.2]{Sno09}), shows that $\rho$ is a member of a strongly compatible system $(\rho_\ell)_\ell$ of representations $\rho_\ell:G_F\to\GL_2(\overline\Q_\ell)$, with $\rho_p=\rho\otimes\overline\Q_p$ and $\rho_\ell|_{G_W}\simeq r_{\pi,\iota_\ell}$ for isomorphisms $\iota_\ell:\overline\Q_\ell\simeq\mathbb C$: for every finite place $v$ of $F$, the Frobenius-semisimple Weil--Deligne representation attached to $\rho_\ell|_{G_{F_v}}$ (via $D_{\mathrm{pst}}$ if $v\mid\ell$) is independent of $\ell$, up to the choice of coefficient embeddings. Each $\rho_\ell$ is finitely ramified and odd, and it is of weight two in the sense of \cite[§1.2]{Sno09}, since its restriction to $G_W$ is $r_{\pi,\iota_\ell}$ with $\pi$ of parallel weight $2$. Since $\pi$ does not have complex multiplication, for all but finitely many $\ell$ the image of $\bar r_{\pi,\iota_\ell}$ contains a conjugate of $\SL_2(\F_\ell)$ \cite{Dim05}. For $\ell\geq5$ the group $\SL_2(\F_\ell)$ is perfect, so it is contained in the image of $G_{W(\zeta_\ell)}$, as $W(\zeta_\ell)/W$ is abelian; hence $\bar\rho_\ell|_{G_{W(\zeta_\ell)}}$, and a fortiori $\bar\rho_\ell|_{G_{F(\zeta_\ell)}}$, is absolutely irreducible.

Choose such a prime $\ell\notin\{2,5,p\}$, different from the residue characteristic of $v_0$ in the second case of (4). Then $\rho_\ell$ satisfies the hypotheses of \cite[Corollary~1.1.2(4)]{Sno09}: either $[F:\Q]$ is odd, or $v_0\nmid\ell$ and $\rho_\ell|_{G_{F_{v_0}}}$ is indecomposable, since its Frobenius-semisimple Weil--Deligne representation is obtained from that of $\rho|_{G_{F_{v_0}}}$ by a change of coefficient embedding. Hence $\rho_\ell$ is the $\lambda_\ell$-adic Tate module of an abelian variety $A/F$ of $\GL_2$-type, for a prime $\lambda_\ell\mid\ell$ of its field of coefficients $K$. The $\lambda$-adic Tate modules of $A$, for the primes $\lambda$ of $K$, form a compatible system whose Frobenius traces agree with those of the system containing $\rho$, through the embedding defined by $\lambda_\ell$. Choosing a $p$-adic embedding of the common trace field compatible with the one defining $\rho$, we get a prime $\lambda\mid p$ of $K$ such that $\rho$ and $\rho_{A,\lambda}$ have the same Frobenius traces, and by the Chebotarev density theorem and the irreducibility of $\rho$, $\rho\otimes\overline\Q_p\simeq\rho_{A,\lambda}\otimes\overline\Q_p$. This concludes the proof.
\end{proof}

%======================================================================

\end{document}